\documentclass[a4paper,12pt]{article}

\usepackage[pdftex]{graphicx}
\usepackage[section] {placeins}
\usepackage{enumerate}
\usepackage[normalem]{ulem}
\usepackage{amsmath}
\usepackage{amsthm}
\usepackage{latexsym}
\usepackage{amsfonts}
\usepackage{amssymb}

\usepackage{tikz-cd}
\usetikzlibrary{graphs}
\usetikzlibrary{arrows.meta}
\usetikzlibrary{shapes.geometric}

\usepackage{authblk}

\usepackage{latexsym,amssymb,amsfonts, amsthm, amsmath}
\usepackage[english]{babel}
\usepackage{bm}
\usepackage{mathrsfs}

\newtheorem{thm}{Theorem}[section]

\newtheorem{lem}[thm]{Lemma}

\theoremstyle{definition} 
 
\newtheorem{exa}[thm]{Example}

\newcommand{\Z}{\mathbb{Z}}

\usepackage{fullpage}

\begin{document}

\title{A Discrete Logarithm Construction for \\ Orthogonal Double Covers of the Complete Graph \\ by Hamiltonian Paths}

\author{M.~A.~Ollis} 

\affil{Marlboro Institute for Liberal Arts and Interdisciplinary Studies, Emerson College, Boston, Massachusetts 02116, USA \\
\texttt{matt\_ollis@emerson.edu}}


\maketitle

\begin{abstract}
During their investigation of power-sequence terraces, 
Anderson and Preece briefly mention a construction of a terrace for the cyclic group~$\Z_n$ when~$n$ is odd and~$2n+1$ is prime; it is built using the discrete logarithm modulo~$2n+1$.  In this short note we see that this terrace  
gives rise to an orthogonal double cover (ODC) for the complete graph~$K_n$ by Hamiltonian paths.  This gives infinitely many new values for which such an ODC is known. 

\bigskip



\end{abstract}

\section{Introduction}\label{sec:intro}

Let~$G$ be a graph with~$n$ vertices.
An {\em orthogonal double cover (ODC)} of the complete graph~$K_n$ by~$G$ is a collection $\mathcal{G} = \{G_i \cong G : 1 \leq i \leq n\}$ that satisfies the following two properties:
\begin{enumerate}
\item {\em Double Cover Property}: Each edge of~$K_n$ lies in exactly two graphs in~$\mathcal{G}$.
\item {\em Orthogonality Property}: Any two distinct graphs in $\mathcal{G}$ have exactly one edge in common.
\end{enumerate}
If $K_n$ has an ODC by~$G$ then~$G$ has exactly~$n-1$ edges.
Figure~\ref{fig:odc9} gives an example of an ODC of~$K_9$ by Hamiltonian paths.
\begin{figure}[tp]
\caption{An orthogonal double cover of~$K_9$ by Hamiltonian paths.}\label{fig:odc9}
$$
\begin{array}{c}
(0,1,4,2,7,5,6,3,8) \\
(1,2,5,3,8,6,7,4,0) \\
(2,3,6,4,0,7,8,5,1) \\
(3,4,7,5,1,8,0,6,2) \\
(4,5,8,6,2,0,1,7,3) \\
(5,6,0,7,3,1,2,8,4) \\
(6,7,1,8,4,2,3,0,5) \\
(7,8,2,0,5,3,4,1,6) \\
(8,0,3,1,6,4,5,2,7)
\end{array}
$$
\end{figure}

Orthogonal double covers have a variety of applications.  See the survey paper~\cite{GGHLL02} for a discussion of these and of the history and development of the problem.  In this note we are concerned with the case when~$G$ is a Hamiltonian path and~$n$ is odd.  This is a case both of independent interest and of application to the almost-Hamiltonian cycle case, one of the original motivators for the study of ODCs.  Theorem~\ref{th:known} captures the current state of knowledge with these constraints.

\begin{thm}\label{th:known}
Let~$n$ be odd with $n>1$.  The complete graph~$K_n$ has an ODC by Hamiltonian paths when it is possible write $n=ab$, where:
\begin{itemize}
\item $a$ has one of the forms $(k^2+1)/2$, $k^2$ or $k^2 + 1$ for some integer~$k$, or 
$$a \in \{3,7,11,15,19, 21, 33, 57, 69, 77, 93  \},$$
\item $b=1$ or $b = q_1q_2\cdots q_m$ for not-necessarily-distinct prime powers~$q_i \equiv 1 \pmod{4}$, with each~$q_i$ either of one of the forms $(k^2+1)/2$, $k^2$ or $k^2 + 1$ for some integer~$k$, or prime with~$q_i < 10^5$.
\end{itemize}
\end{thm}

\begin{proof}[Proof note]
Given an ODC of~$K_a$ by Hamiltonian paths (or the trivial case~$a=1$), we may combine this with a ``2-colorable" ODC for $K_{q}$ by Hamiltonian paths, where~$q$ is a prime power congruent to~$1 \pmod{4}$, to obtain an ODC of~$K_{aq}$ by Hamiltonian paths.  In the statement of the theorem, the values of~$q_i$ are those for which such a 2-colorable ODC is known.  There is an ODC of~$K_a$ by Hamiltonian paths for each of the non-trivial values of~$a$ and, in conjunction with the possibilities for~$b$, the statement covers all odd values~$n$ for which ODCs of~$K_n$ by Hamiltonian paths are known.

For the full definitions and constructions that prove this result, see~\cite{GGHLL02,HLL08} and the references therein.
\end{proof}

A ``terrace", to be defined in the next section for the cases we need, is an arrangement of the elements of a group such that adjacent elements have particular difference/quotient properties.  Terraces for arbitrary groups were formally introduced by Bailey in~1984~\cite{Bailey84}. However, the concept for cyclic groups (which is all we require in this note) is implicit in work of Williams in 1949~\cite{Williams49} and, arguably, in that of Lucas and Walecki in 1892, see~\cite{Alspach08}.  

A terrace for a group of order~$n$ gives rise to collection of~$n$ Hamiltonian paths in the complete graph~$K_n$ that has the double cover property (equivalently, a {\em row quasi-complete latin square of order~$n$})~\cite{Bailey84}.   A terrace may have additional properties that mean we in fact have an orthogonal double cover of~$K_n$ by Hamiltonian paths~\cite{AL91,HN91,Leck98}.  These additional properties are rare and it seems to be a difficult problem to find terraces that have them.

 In~\cite{AP03a}, the first in a series of papers investigating ``power-sequence terraces" for cyclic groups (see~\cite{AOP17} for references for the full series), Anderson and Preece give a construction for a terrace of~$\Z_n$ using discrete logarithms modulo~$2n+1$ that applies when~$n$ is odd and~$2n+1$ is prime.  In the next section we prove that this terrace has the additional properties required to give an ODC of~$K_n$ by Hamiltonian paths.

We are thus able to augment the list of values of~$n$ for which~$K_n$ is known to have an ODC by Hamiltonian paths as follows.

\begin{thm}\label{th:main}
If~$n$ is odd and~$2n+1$ is prime, then~$K_n$ has an orthogonal double cover by Hamiltonian paths.
\end{thm}

Any value of~$n$ that satisfies the conditions of Theorem~\ref{th:main} may be added as a possibility for~$a$ in Theorem~\ref{th:known}.  Each new value that is congruent to~$3 \pmod{4}$ gives rise to infinitely many new values for which an ODC of the complete graph by Hamiltonian paths is known.  These are straightforward to describe: for each prime~$p$ with $p \equiv 7 \pmod{8}$ we obtain the value~$(p-1)/2$, which is a new value when~$p \geq 47$ and $(p-1)/2$ is not of one of the forms $(k^2+1)/2$,~$k^2$ or $k^2 + 1$ for some integer~$k$.

Given the possibilities for~$b$ in Theorem~\ref{th:known}, values given by Theorem~\ref{th:main} that are congruent to~$1 \pmod{4}$ do not immediately cover as many new cases. 

If the value of~$n$ is prime (in which case it is a {\em Sophie Germain prime}) and it is not of one of the forms $(k^2+1)/2$ or $k^2 + 1$ for some integer~$k$, then it gives infinitely many new values with~$n \equiv 1 \pmod{4}$ (provided $n > 10^5$) and $n \equiv 3 \pmod{4}$ (provided~$n > 20$).  Sophie Germain primes are well-studied and many such are known.  See, for example,~\cite[A005384]{OEIS} and the linked files. 

We also get values congruent to~$1 \pmod{4}$, each corresponding to a infinite family, by taking $n$ to be the product of an even number of distinct primes, each congruent to~$3 \pmod{4}$ (subject, of course, to~$2n+1$ being prime).  These are new when~$n > 100$ and is not of the form $(k^2+1)/2$ or $k^2 + 1$ for some integer~$k$.

\section{The Construction}\label{sec:const}

We shall work in $\Z_n = \{0, 1, \ldots, n-1\}$, the additively-written cyclic group of order~$n$. In this section,~$n$ is always odd; write~$n=2m+1$.  
When we consider the complete graph~$K_n$, it is convenient to consider it as having its vertices labeled with the elements of~$\Z_n$.  

Define the {\em length} of the edge between vertices~$x$ and~$y$ to be $\{y-x, x-y\}$.  For brevity, we sometimes abbreviate this to~$\pm (y-x)$.  Let~$\mathbf{h}$ be a Hamiltonian path in~$K_n$.  If each edge length~$\pm \ell$ for $\ell \in \{1, 2, \ldots, m\}$ occurs exactly twice in~$\mathbf{h}$ then~$\mathbf{h}$ is a {\em terrace} for~$\Z_n$~\cite{Bailey84}. 

Define the {\em distance} between two edges $e_1 = (x_1,y_1)$ and $e_2 = (x_2,y_2)$ of the same length to be~$\pm k$ if either 
$$\{x_1 + k, y_1 + k\} = \{x_2, y_2\} \textrm{ \ or \ } \{x_1 - k, y_1 - k\} = \{x_2, y_2\}.$$
If a terrace~$\mathbf{h}$ has the property that every possible distance~$\pm k$, for $k \in \{1, 2, \ldots, m\}$, appears between the pairs of edges of the same length, then $\mathbf{h}$ is an {\em ODC-starter} for $\Z_n$~\cite{HN91,Leck98}.

\begin{lem}\label{lem:starter}{\rm \cite{HN91,Leck98}}
If $\mathbf{h}$ is an ODC-starter in~$K_n$, then the translates of~$\mathbf{h}$ form an orthogonal double cover of~$K_n$ by Hamiltonian paths.
\end{lem}

\begin{exa}\label{ex:starter}
Let~$n=9$ and consider the Hamiltonian path~$\mathbf{h} = (0,1,4,2,7,5,6,3,8)$ in~$K_9$.  This is a terrace; the sequence of edge lengths is
$$(\pm 1, \pm 3, \pm 2, \pm 4, \pm 2, \pm 1, \pm 3, \pm 4).$$
Further, it is a ODC-starter.  The distances associated with the edge lengths~$\pm 1$, $\pm 2$, $\pm 3$ and $\pm 4$ are $\pm 4$, $\pm 3$, $\pm 2$ and  $\pm 1$ respectively.
Hence, by Lemma~\ref{lem:starter}, the translates of~$\mathbf{h}$ are an orthogonal double cover of~$K_9$ by Hamiltonian paths; this is the ODC displayed in Figure~\ref{fig:odc9}.  
\end{exa}

Here is Anderson and Preece's complete discussion of their construction in~\cite{AP03a} (citation updated to align with the numbering of the present paper):

\vspace{2mm}
\setlength{\leftskip}{1cm}
\setlength{\rightskip}{1cm}
\begin{small}

...there is also a general Galois field construction for~$\Z_n$ terraces for any odd integer~$n$ such that~$2n+1$ is prime; some of the mathematics for this is in~\cite{PVHRV95}.  Let~$x$ be a primitive root of $\Z_{2n+1}$.  Working in $\Z_{2n+1}$, obtain the quantities $c_i$ ($i = 1, 2, \ldots, n$) defined by $x^{c_i} = i$.  Now reduce $c_i$ modulo~$n$, to give  $d_i$ ($i = 1, 2, \ldots, n$).  Then $(d_1, d_2, \ldots, d_n)$ is a $\Z_n$ terrace.  With $n=5$ and $x=2$ we obtain the $\Z_5$ terrace $(0,1,3,2,4)$, and with $n=9$ and $x=2$ we obtain the $\Z_9$ terrace $(0,1,4,2,7,5,6,3,8)$.

\end{small}
\setlength{\leftskip}{0pt}
\setlength{\rightskip}{0pt}
\vspace{2mm}
We shall prove Theorem~\ref{th:main} by showing that not only is Anderson and Preece's discrete logarithm construction a terrace (Theorem~\ref{th:ap_terrace}), it is also an ODC-starter (Theorem~\ref{th:ap_odc}).  Notice that Example~\ref{ex:starter} verifies this for the case given with~$n=9$ and~$x=2$. 

Here is some notation, which will hold for the remainder of the paper.  Let~$g$ be a primitive root of~$\Z_{2n+1}$, where~$n$ is odd.  We take an element~$x \in \Z_{2n+1}$ to~$\widehat{ \log_g(x) } \in  \Z_n$ in two steps, where the first is the discrete logarithm $\log_g$ with respect~$g$ and lies in~$\Z_{2n}$ and the second is the natural projection~$\Z_{2n} \rightarrow \Z_n$, which we denote in general by~$y \mapsto \widehat{y}$.
Then Anderson and Preece's construction is given by~${\bf d}_g = (d_1, d_2, \ldots, d_n)$, where $d_i =  \widehat{ \log_g(i) }$ for each~$i$.

We need one more concept before the results.  Let ${\bf a} = (a_1, a_2, \ldots, a_{2n})$ be an arrangement of the elements of~$\Z_{2n}$ and let~$b_i = a_{i+1} - a_i$ for~$i$ in the range~$1 \leq i < 2n$.  If the elements~$b_i$ are all of the non-zero elements of~$\Z_n$ and~$b_i = -b_{2n-i}$ for~$i$ in the range~$1 \leq i < n-1$ then~${\bf a}$ is a {\em symmetric directed terrace} for~$\Z_{2n}$ and the sequence $(b_1, b_2, \ldots, b_{2n-1})$ is a {\em symmetric sequencing}.  In~\cite{Anderson76} it is shown that if $(a_1, a_2, \ldots, a_{2n})$ is a symmetric directed terrace then $(\widehat{a_1}, \widehat{a_2}, \ldots, \widehat{a_{n}})$ is a terrace.  (This concept, and result, works more generally in groups that have exactly one involution, and it is always possible to go from the terrace to a symmetric sequencing; see, for example,~\cite{OllisSurvey}.)

\begin{thm}\label{th:ap_terrace}
Anderson and Preece's discrete logarithm construction gives a terrace for each odd~$n$ such that~$2n+1$ is prime.
\end{thm}

\begin{proof}
It is sufficient to show that the sequence
${\bf c}_g = (c_1, c_2, \ldots, c_{2n})$ given by $c_i = \log_g(i)$ is a symmetric directed terrace for~$\Z_{2n}$.  Let~$b_i = c_{i+1} - c_i$ for~$i$ in the range~$1 \leq i < 2n$.

As~$g$ is a primitive root, when~$1 \leq i, j \leq 2n$ with~$i \neq j$ we have 
$$c_i = \log_g(i) \neq \log_g(j) = c_j$$
and so~${\bf c}_g$ is an arrangement of the elements of~$\Z_{2n}$.

The difference~$b_i$ is given by
$$b_i = \log_g(i+1) - \log_g(i) = \log_g \left( \frac{i+1}{i} \right).$$ 
When~$1 \leq i, j \leq 2n$ with~$i \neq j$, we have $(i+1)/i \neq (j+1)/j$ in~$\Z_{2n+1}$ and hence $b_i \neq b_j$ in~$\Z_{2n}$.

Let~$1 \leq i < n$.  In~$\Z_{2n+1}$ we have
$$ \frac{2n-i+1}{2n-i} = \frac{i}{i+1} =  \left( \frac{i+1}{i} \right)^{-1}.$$
Therefore, in~$\Z_{2n}$,
$$b_{2n-i} = \log_g \left( \frac{2n-i+1}{2n-i} \right)
 = \log_g \left( \left( \frac{i+1}{i} \right)^{-1} \right)
= - \log_g \left( \frac{i+1}{i} \right) = -b_i.
$$

Hence~${\bf c}_g$ is a symmetric directed terrace for~$\Z_{2n}$ and~${\bf d}_g$ is a terrace for~$\Z_n$.
\end{proof}

\begin{thm}\label{th:ap_odc}
Anderson and Preece's discrete logarithm construction gives an ODC-starter for each odd~$n$ such that~$2n+1$ is prime.
\end{thm}

\begin{proof}
Theorem~\ref{th:ap_terrace} shows that Anderson and Preece's discrete logarithm construction ${\bf d}_g$ is a terrace.  It remains to show that between the pairs of edges of the same length we realise all possible distances.  Recall that~$n = 2m+1$.

Fix~$k$ in the range~$1 \leq k \leq m$.  We want to show that there is a pair of edges of the same length that are at distance~$\pm k$.

Let~$x = g^k$ in~$\Z_{2n+1}$ and let $u = (1-x)/(1+x)$, which is well-defined as $x = -1$ only when~$k = n = 2m+1$.   
Let $i = (u - 1)^{-1}$ and $j = xi$ in~$\Z_{2n+1}$.  The pairs $(i, i+1)$ and $(j, j+1)$ in~$\Z_{2n+1}$ become the edges
$$e_i = \{ \widehat{\log_g(i)}, \widehat{\log_g(i+1)} \} \text{ \ \ and \ \ } e_j = \{ \widehat{\log_g(j)}, \widehat{\log_g(j+1)} \}$$
in the terrace for~$\Z_n$.  We claim that these edges have the same length and are at distance~$\pm k$.

As $i = (u - 1)^{-1}$, we have $i^{-1} = u-1$ and can evaluate
$$\frac{i+1}{i} = 1 + \frac{1}{i} = 1 + (u-1) = u$$
in~$\Z_{2n+1}$.
As~$j=xi$ and writing~$u-1$ in terms of~$x$ as $-2x/(1+x)$ we also have
$$ \frac{j+1}{j} = 1 + \frac{1}{xi} = 1 + \frac{u-1}{x}
= \frac{x-1}{1+x} = -u$$
in~$\Z_{2n+1}$.

In~$\Z_{2n+1}$, we have $g^n = -1$ and so in~$\Z_n$ we have $\widehat{\log_g(y)} = \widehat{\log_g(-y)}$ for all~$y$.  In particular, the lengths of~$e_i$ and~$e_j$ are both $\widehat{\log_g(u)} = \widehat{\log_g(-u)} $.

We now show that~$e_1$ and~$e_2$ are at distance~$\pm k$.
First, since $j = xi = g^ki$, we have
$$\widehat{\log_g(i)} + k = \widehat{\log_g(j)}$$
in~$\Z_n$.  Second, note that $j+1 = -x(i+1)$ in~$\Z_{2n+1}$, hence
$$\log_g(j+1) = \log_g(-1) + k + \log(i+1) = n +  k + \log(i+1)$$
in~$\Z_{2n}$.  This gives
$$ \widehat{\log_g(j+1)} = \widehat{ \log_g(i+1)} + k$$
in~$\Z_n$, and so the edges are indeed at distance~$\pm k$.

Therefore, Anderson and Preece's discrete logarithm construction has all the properties required to be an ODC-starter.
\end{proof}

An orthogonal double cover of~$K_n$ by Hamiltonian paths, where~$n$ is odd and~$2n+1$ is prime, now exists by Lemma~\ref{lem:starter}.  This proves Theorem~\ref{th:main}.

\begin{exa}
Let~$n=15$.  Then~$2n+1 = 31$, which is prime and has~3 as a primitive root.  Using this, the Anderson and Preece discrete logatirhm construction gives the
following ODC-starter in~$\Z_{15}$:
$$(0,9,1,3,5,10,13,12,2,14,8,4,11,7,6).$$
Its sequence of edge lengths is
$$(\pm6,\pm7,\pm2,\pm2,\pm5,\pm3,\pm1,\pm5,\pm3,\pm6,\pm4,\pm7\pm4,\pm1).$$
The distances between the pairs of edges of lengths $\pm1, \pm2, \ldots, \pm7$ are 
$$ \pm6, \pm2, \pm4, \pm3, \pm7, \pm1, \pm5$$  
respectively.
\end{exa}

\end{document}